\documentclass[runningheads]{PMS2026}

\usepackage[ansinew]{inputenc}
\usepackage[T1]{fontenc}
\usepackage[top=3cm,bottom=3.5cm,left=4cm,right=3.8cm,dvips]{geometry}
\usepackage{nopageno}

\usepackage{index}
\usepackage{fancyhdr}

\usepackage{indentfirst}
\usepackage{sectsty}
\allsectionsfont{\usefont{OT1}{phv}{bc}{n}\selectfont \normalsize \bf }
\usepackage[dvips]{graphicx}
\usepackage{subfigure}
\usepackage{url}
\usepackage{harvard}
\usepackage{amssymb}
\usepackage{epsfig}

\usepackage{color}

\usepackage{amsmath}
\usepackage{anyfontsize}

\renewcommand{\figurename}{Model}

\makeindex

\begin{document}
\pagenumbering{gobble}
%\pagenumbering{roman}
\pagestyle{headings}
\mainmatter

%\tableofcontents

\title{Scheduling of Star Observations under Uncertain Conditions: A Comparison of Models and Solvers}

\author{Thomas Rahab Lacroix, Pierre Lemaire, Nadia Brauner}

\institute{Univ. Grenoble Alpes, CNRS, Grenoble INP, G-SCOP, 38000 Grenoble, France\\
\email{\{thomas.rahab-lacroix, pierre.lemaire, nadia.brauner\}@grenoble-inp.fr}
}

\maketitle
% Please, write aprofit your names here under, as you want it to appear in the author's index.
\index{Rahab Lacroix Thomas}\index{Lemaire Pierre}\index{Brauner Nadia}

\noindent
\textbf{Keywords:} scheduling, optimization under uncertainty, robustness, min-max regret, arc flow, solver comparison.

\section {Scheduling with uncertainty on the number of machines}

We consider a scheduling problem for star observations using a telescope: $M$ identical nights are available to observe $J$ stars that have been identified as of interest by astronomers. Each star $j \in \{1,2,...,J\}$ yields a profit $w_j$ if it is observed; in that case, the observation must respect a visibility window $[r_j,d_j]$ and a minimum observation duration of at least $p_j$. The objective is to maximize the total profit of the observations that are actually performed. By interpreting nights as machines and stars as jobs to be processed, this problem is exactly $PM|r_j|\sum w_j U_j$ in Graham's notation, a classic scheduling setting of parallel machines.

In the practical star-observation problem, the ability to observe depends on meteorological and atmospheric conditions, which are not known when the schedule is computed. A possible modeling is to consider that each night can be either perfect and everything is observable, either terrible and nothing is observable. The challenge is therefore to propose a "high-quality" schedule without knowing the actual number of nights (machines) that will ultimately be available \cite{stein_zhong}; this contrasts with classical scheduling models, where uncertainty typically affects the job-related data.

This problem can be viewed from the perspective of stochastic optimization when probabilities over the number of available nights are known \cite{moi}. Since such probabilities are difficult to obtain in practice, we focus here instead on evaluating the worst case within the robust optimization paradigm.

\section{{\em Min-max regret} formulation}

As the nights are identical, when one appears observable, a schedule must execute its night with the highest profit. So we can assume that all schedules have their nights sorted by decreasing profit. We denote by $\mathcal{X}$ the set of feasible schedules over $M$ nights. The objective is therefore to propose a schedule $x \in \mathcal{X}$ computed for $M$ nights such that for each $m = 1,...,M$, its restriction to its first $m$ nights is not too far from the optimal schedule over $m$ nights. Let $x(m)$ denotes the cumulative profit of schedule $x$ over its first $m$ nights, and let $OPT_m \in \mathcal{X}$ be an optimal schedule over $m$ nights for the problem without uncertainty. The difference between these two values is the regret, and the minimization of the maximum regret \cite{robust_scheduling} can then be written as (\ref{1easy}). (\ref{1easy}) can be rewritten as a decision problem (\ref{eq:prbldecision}) by introducing a bound $R$:

\begin{equation}\label{1easy}
    \underset{x\in \mathcal{X}}{\min} \quad \Biggl( \underset{m \in \{1, 2..., M\}}{\max} \quad \Bigl( OPT_m(m)  - x(m) \Bigr) \Biggr)
\end{equation}

\begin{equation}\label{eq:prbldecision}
    \exists x   \in \mathcal{X}, \quad \forall y_1, y_2..., y_M \in \mathcal{X}, \quad \forall m\in  \{1, 2..., M\}, \quad y_m(m) - x(m) \leq R
\end{equation}

Indeed, for any schedule $y\in \mathcal{X}$, we have $OPT_m(m)\geq y(m)$, and there exists one for which equality holds. Therefore, the above formulation corresponds precisely to the decision version of the maximum regret minimization problem. The alternation of an ``exists'' and a  ``for all'' quantifier places this problem in the complexity class $\Sigma_2^p$ \cite{sigma_2p}. Moreover, the quantifiers can be inverted, so the problem is also in $\Pi_2^P$. Besides, the problem can be formulated as a MILP with an objective function different than just the regret, and thus belong to $\Delta_2^P$. However, the exact complexity is not known yet.

\section{Models for the computation of an optimal schedule}
To compute the regret, it is essential to compute as efficiently as possible an optimal schedule over $m$ machines, for each $m=1,..., M$. We consider three possible Integer Linear Programming (ILP) formulations for this scheduling problem: Start Time (ST), Time Discretization (TD), and Arc Flow (AF) and compare their performances.

\bigskip

The classical Start Time formulation \cite{start_time} involves $O(J^2)$ binary variables (see Model \ref{fig:start_time}). Its simplicity is balanced by poor solving capacity due to the use of big-$M$ (here a $T = \max_j(d_j)$) constraints and the quickly increasing number of binary variables.

\medskip

In the Time Discretization formulation \cite{time_discr}, a night is discretized into $T+1$ instants. This model involves $O(MTJ)$ binary variables (see Model \ref{fig:time_discretization}). Its flexibility and ability to solve quite large instances is limited when the number of nights, or the number of instants per night, is too high.

\medskip

The Arc Flow formulation \cite{Kramer} models the scheduling problem as a flow in a graph. Time is discretized into $T+1$ instants, and the vertices correspond to the time instants from 0 to $T$. For each star $j$ and each time $t\in[r_j,d_j-p_j]$, an arc from $t$ to $t+p_j$ with weight $w_j$ is created. Additional arcs are used to represent idle times. A path from vertex 0 to vertex $T$ in this graph corresponds to the schedule of a single night. The objective is to find a maximum-weight flow from vertex 0 to vertex $T$ that satisfies additional constraints: each star can be scheduled at most once, and the total flow leaving 0 is bounded by $M$.

The Arc Flow model involves $O(TJ)$ binary variables (see Model \ref{fig:arc_flow}) and is considered in the literature very powerful. However this formulation cannot be used to solve other more general variants of the problem, e.g. when nights are not identical. Because there is no information about the night when an observation is scheduled, this formulation cannot be used when partial evaluations are required (e.g. robust or stochastic optimization).

\begin{figure}[htbp]
\footnotesize
$a_{nj} = 1$ if observation $j$ is scheduled night $n$, and $0$ otherwise

$b_{ij} = 1$ if observation $i$ is scheduled before observation $j$ the same night, and $0$ otherwise

$c_j = $ start time of observation $j$
\[
\begin{aligned}
\text{Maximize} \quad     & \sum_{n=1}^M \sum_{j=1}^J a_{nj} w_j     & & \\%PL: \textbf{Start Time formulation } (3)\\
\forall j, \quad          & \sum_{n=1}^M a_{nj} \leq 1               & & \text{observation $j$ is scheduled at most once} \\
\forall j, \quad          & r_j \leq c_j \leq d_j - p_j              & & \text{observation $j$ is scheduled in its time window} \\
\forall i \neq j, \quad   & c_i + p_i \leq c_j + T(1-b_{ij})         & & \text{ensure the no overlap of observations $i$ and $j$} \\
\forall i \neq j,n, \quad & b_{ij} + b_{ji} \geq a_{nj} + a_{ni} - 1 & & \text{enforce an order for observations $i$ and $j$} \\
\forall i,j,n, \quad      & a_{nj} \in \{0,1\}, \indent b_{ij} \in \{0,1\},  & &  c_j \geq 0
\end{aligned}
\]
\caption{Start Time (Wagner 1959) }\label{fig:start_time}
\end{figure}

\begin{figure}[htbp]
\footnotesize
$ s_{ntj} = 1 $ if observation $j$ starts at instant $t$ the night $n$, and $0$ otherwise
\[
\begin{aligned}
\text{Maximize} \quad & \sum_{n=1}^M \sum_{j=1}^J w_j \sum_{t=0}^T s_{ntj} \\
\forall j, \quad      & \sum_{n=1}^M \sum_{t=0}^T s_{ntj} \leq 1           & & \text{observation $j$ is scheduled at most once} \\
\forall j,n, \quad    & \forall t \in [\![0,r_j -1]\!], s_{ntj} = 0        & & \text{observation $j$ starts after $r_j$} \\
\forall j,n, \quad    & \forall t \in [\![d_j - p_j +1, T]\!], s_{ntj} = 0 & & \text{observation $j$ ends before $d_j$} \\
\forall t,n, \quad    & \sum_{j=1}^J \sum_{u=t-p_j+1}^t s_{nuj} \leq 1     & & \text{at most one observation an instant $t$} \\
\forall j,t,n \quad   & s_{ntj} \in \{0,1\}
\end{aligned}
\]
\caption{Time Discretization (Sousa et. al. 1992) }\label{fig:time_discretization}
\end{figure}

\begin{figure}[htbp]
\footnotesize
$ f_{tj} = 1 $ if observation $j$ starts at instant $t$, and $0$ otherwise
\[
\begin{aligned}
\text{Maximize} \quad   & \sum_{j=1}^J \sum_{t=0}^T f_{tj} w_j \\
\forall j \neq 0, \quad & \sum_{t=0}^T f_{tj} \leq 1                                        & & \text{observation $j$ is scheduled at most once} \\
\forall j, \quad        & \forall t \in [\![0,r_j -1]\!], f_{tj} = 0                        & & \text{observation $j$ starts after $r_j$} \\
\forall j, \quad        & \forall t \in [\![d_j - p_j +1, T]\!], f_{tj} = 0                 & & \text{observation $j$ ends before $d_j$} \\
\forall t, \quad        & \sum_{j=0}^J f_{tj} - f_{(t-p_j)j} = \begin{cases}
                                                                    M \text{ if } t = 0\\
                                                                    -M \text{ if } t=T\\
                                                                    0 \text{ else}
                                                                \end{cases}                 & & \text{flow conservation} \\
\forall j,t, \quad      & f_{tj} \in \{0,1\} \text{ except } f_{t0} \in \mathbb{R}^+
\end{aligned}
\]
\caption{Arc flow (Kramer et. al. 2020)}\label{fig:arc_flow}
\end{figure}

\section{Comparison of solvers and models}

Astronomers provided us with a list of 800 stars with their profits. Then, 120 dates were selected uniformly over the 10 last years and time windows and processing times for all stars were computed for each situation. As consecutive nights are quite similar in terms of time windows, an instance over $M$ nights is built by repeating $M$ times the same situation (with $M$ taking its value in $\{1, 3, 7, 14, 30, 100\}$). For each instance, we keep all visible stars (about 400 to 600 among the 800 stars). At the end, we obtain 720 distinct realistic instances.

Those different instances are solved multiple times. With the 3 models presented before: ST, TD and AF but also with 2 solvers: HiGHS and CPLEX \cite{solver}. This brings it to a total of 4320 resolutions and each of them is allowed 2 minutes of computation on a server equipped with an Intel(R) Xeon(R) Silver 4210R CPU @ 2.40GHz, 20 physical cores, and 40 threads.

\begin{table}[htbp]
\centering
\hspace*{-30pt}
    \begin{tabular}{c|cccccc|cccccc|}
           & \multicolumn{6}{c|}{Number of instances where} 
           & \multicolumn{6}{c|}{Mean time (s) to prove} \\
           & \multicolumn{6}{c|}{optimality is proven over 120 dates} 
           & \multicolumn{6}{c|}{optimality when proven under 2min} \\
           & 1       & 3       & 7       & 14      & 30      & 100    & 1        & 3          & 7        & 14       & 30       & 100 \\
        \hline
        ST &  0/0    & 0/0     &  0/0    & 0/0     & 0/0     &  0/0   &  \_/\_   & \_/\_      & \_/\_    &  \_/\_   & \_/\_    & \_/\_  \\
        TD & 99/120  & 2/120   &  0/3    & 0/0     & 0/0     &  0/0   & 64.7/7.6 & 122.9/30.4 & \_/108.7 & \_/\_    & \_/\_    & \_/\_   \\
        AF & 118/120 & 120/119 & 120/120 & 120/120 & 118/120 & 92/120 & 43.6/4.3 & 16.5/3.8   & 20.7/4.3 & 30.5/4.5 & 38.7/5.0 & 47.0/6.0
    \end{tabular}
\smallskip
\caption{Performances of the different models and solvers. Each cell contains two values that are obtained for the above metric with the alongside formulation for the above number of nights with the solvers HiGHS/CPLEX.} \label{table:metric_model}
\end{table}

A clear results from Table~\ref{table:metric_model} about the models is that the AF formulation is by far the best tool to solve the $OPT_m$ problem. Moreover, with the actual characteristics of our instances, the TD formulation is strictly better than the ST formulation.

About the solvers, Table~\ref{table:metric_model} shows that CPLEX is faster than HiGHS to prove optimality no matter the formulation. In average, CPLEX is 6 times faster than HiGHS. This ratio is low considering that CPLEX is a commercial solver and HiGHS an open source solver.

\section{Conclusion}

We presented three models to solve the $PM|r_j|\sum w_j U_j$ problem. Their performances are compared on realistically sized instances (up to 100 nights and 600 stars) using different solvers that are also compared. We conclude that to solve this problem, the Arc Flow formulation is by far the most efficient. In the perspective of solving a robust or stochastic formulation, the Time Discretization model seems more appropriate than the Start Time model. Besides, reducing the size of an instance by using only stars selected by Arc Flow for an optimal resolution is a promising heuristic.

To conclude, the next step is to add metrics to analyse the performances of the models and solvers, such as the number of instances with a gap to optimality lower than 5\% and the number of instances reaching the best known solution obtained with the same model. Moreover, new commercial and open source solvers will be add to the comparison. Experiments with CBC, Gurobi, and GLPK are currently in progress.

%\printindex
\end{document}